\documentclass[11pt,a4paper, twoside]{article}
\usepackage[latin9]{inputenc}
\usepackage{amsmath,amssymb,amsthm,amscd,hyperref}
\usepackage[]{fontenc}
\usepackage{xy}
\usepackage{enumerate}
\usepackage{graphicx}

\newif\ifpdf\ifx\pdfoutput\undefined\pdffalse\else\pdfoutput=1\pdftrue\fi

\xyoption{all} \numberwithin{equation}{section}
\newtheorem{theorem}[equation]{Theorem}

\newtheorem{proposition}[equation]{Proposition}
\newtheorem{lemma}[equation]{Lemma}

\theoremstyle{definition}

\newtheorem{definition}[equation]{Definition}

\newtheorem{rem}[equation]{Remark}

\newcommand{\nc}{\newcommand}
\nc{\cH}{\mathcal{H}} \nc{\cG}{\mathcal{G}}
\nc{\cC}{\mathcal{C}}
\nc{\cO}{\mathcal{O}}
\nc{\cI}{\mathcal{I}}
\nc{\cA}{\mathcal{A}}
\nc{\cB}{\mathcal{B}} \nc{\cY}{\mathcal{Y}} \nc{\cK}{\mathcal{K}}
\nc{\cX}{\mathcal{X}} \nc{\cS}{\mathcal{S}} \nc{\cE}{\mathcal{E}}
\nc{\cF}{\mathcal{F}} \nc{\cZ}{\mathcal{Z}} \nc{\cQ}{\mathcal{Q}}
\nc{\cN}{\mathcal{N}} \nc{\cP}{\mathcal{P}} \nc{\cL}{\mathcal{L}}
\nc{\cM}{\mathcal{M}} \nc{\cR}{\mathcal{R}} \nc{\cT}{\mathcal{T}}
\nc{\cW}{\mathcal{W}} \nc{\cU}{\mathcal{U}} \nc{\cD}{\mathcal{D}}
\nc{\cJ}{\mathcal{J}} \nc{\cV}{\mathcal{V}}
\nc{\fr}{{\rightarrow}}
\nc{\rd}{red.deg}

\newcommand{\pr}{\mathbb P}
\newcommand{\F}{\mathbb F}

\newcommand{\rk}{\mbox{\upshape{rank}}}

\newcommand{\kapp}{\mathbb C}

\title{Locally non-trivial fibred surfaces with maximal unitary rank}
\author{Lidia Stoppino}

\begin{document}
\maketitle
\begin{abstract}Let $f\colon S\to B$ a locally non-trivial fibred surface with fibres of genus $g$. 
Let $u_f$ be its unitary rank, i.e. the rank of the flat unitary part in the second Fujita decomposition. 
We study in detail the case when $u_f$ is maximal,  i.e. $u_f=g-1$. 
In this case necessarily $g\leq 6$, but examples in genus $5$ and $6$ are not known, and conjecturally do not exist.
We prove a strong slope inequality for these extremal cases.
We then use this inequality, together with results on trigonal curves, to give new constraints on the case $g=6$, $u_f=5$.
In particular, we prove that the index of the surface is always strictly positive and give strong limitations on the possible classes of the relative canonical divisor.
\end{abstract}

{Fibred varieties, slope inequality, families of Jacobians, Hodge bundle, unitary rank, relative irregularity}
[MSC Classification]{14J10, 14D06, 14D07, 14G35, 14H40}

\section{Introduction}
Let $f\colon S\to B$ be a relatively minimal locally non-trivial fibred surface of genus $g\geq 2$.
Consider two invariants closely related to the two Fujita decompositions of the Hodge bundle (see Theorem \ref{teoremone2}): (1) the relative irregularity $q_f:= q(S)-g(B)$; (2) the unitary rank $u_f$, i.e. the rank of the unitary summand in the second Fujita decomposition.
We have that $q_f\leq u_f$, and the inequality can be strict, also up to base changes (see Remark \ref{rem: confronto}). 
In \cite{Beau}, Beauville proved that $q_f\leq g$ and equality holds if and only if the fibration is birationally trivial. 
It can be  easily seen (see Remark \ref{rem: ossqf}) that also $u_f\leq g$, and equality holds if and only if $f$ is locally trivial. 
So, for a locally non-trivial fibration the maximal possible value for $u_f$ is $g-1$. 
From a result of the author in \cite{Sto24}, we know that in this case $g$ has to be smaller or equal to $6$. 

There are several examples in genus $2\leq g\leq 4$:  the first one, in genus $4$, is due to Pirola \cite{pirola}, then a complete classification for $g=2,3$ appeared  in  \cite{moller}, and other examples in genus $4$ in \cite{AP}.
No examples in genus $5$ or $6$ are known, the general expectation being that they do not exist.
 
Thanks to  the work of many authors,  several conditions have been found for locally non-trivial fibrations with $4\leq g\leq 6$ and  $u_f=g-1$: see Theorem \ref{teo: estremale}.
In particular these fibrations  always are trigonal of special Maroni invariant, and  never are Kodaira fibrations.

\medskip

In this note, we firstly prove an extremely strong slope inequality for these fibrations. 
Let $f\colon S\to B$ a relatively minimal locally non-trivial fibred surface of genus $g\geq 2$ such that $u_f= g-1$. 
Let $\cA$ be the ample line bundle in the second Fujita decomposition of the Hodge bundle $f_*\omega_f$ \eqref{second fujita}, where $\omega_f=\omega_S\otimes f^*\omega_B^{-1}$ is the relative canonical line bundle.
Consider the sheaf homomorphism $\eta$
$$
\xymatrix{f^*\cA\ar@/^1pc/[rr]^\eta \ar[r]& f^*f_*\omega_f\ar[r]&\omega_f,}
$$
where the second homomorphism is the evaluation homomorphism.  
Let  $Z$ be the divisorial (schematic) base locus  associated to $\eta$.
We call  $F$ a general smooth fibre of $f$.
Let  $K_f=K_S-f^*K_B$ be the relative canonical divisor, so that $\omega_f\cong \cO_S(K_f)$.

\smallskip

\begin{proposition}(Proposition \ref{prop: bound estremale})
With the above notations,
let $Z=\sum^k_{j=1}m_jG_j+Z'$, where the $G_j$'s are the irreducible components of $Z$ which are not vertical with respect to $f$, and the $m_j$'s are their respective multiplicities. 
Let $b_j:=G_j\cdot F$. 

We have the following inequality:
\begin{equation}\label{eq: meglio}
K_f^2\geq \left(2(g-1)+ \sum_{j=1}^k\frac{m_j}{m_j+1}b_j\right)\chi_f.
\end{equation}
\end{proposition}

 This result is achieved by making a careful use of the special  form that the relative canonical divisor assumes for the extremal case. 
 Indeed,  in these cases (see Lemma  \ref{lem: serve})
 $$K_f\sim  Z+f^*\cA\equiv Z+\mu_1F,$$ where $\sim $ denotes linear equivalence and $\equiv $ numerical equivalence.
In order to  prove this slope inequality we do not use Xiao's method, but rather an ad-hoc computation that crucially uses the Grothendieck-Riemann-Roch Theorem and some arguments in line with \cite[Sec.2]{konno}. 
The multiplicities $m_j$ of the horizontal components of $Z$ play an important r\^ole in this result: the smaller they are the better the bound becomes, as we see in Remark \ref{rem: migliore}.

Then we use this strong slope inequality, combined with the inequality $K_f^2<12 \chi_f$ coming from point (3) of Theorem \ref{teo: estremale}, to exclude many possibilities for $Z$, and hence for the  class of $K_f$. In particular we obtain the following:

\smallskip

\begin{theorem} [Proposition \ref{prop: 6}, Proposition \ref{prop: elenco}, Proposition \ref{prop: ram}]
Let $f\colon S\to B$ a relatively minimal  locally non-trivial  fibred surface of genus $6$ and with $u_f=5$. Then:
\begin{enumerate} 
\item[(1)] the surface $S$ has positive index $\tau(S)>0$.
\item[(2)] The divisor $Z$ has a horizontal component $G$ of multiplicity at least 5. 
\item[(3)] There is a list of the possible isomorphism class of $K_f$, given in in Proposition \ref{prop: elenco}.
\item[(4)] If the multiplicity of $G$ is at least $6$ then $G$ is a section of $f$ and  the point $G\cap F$ on the general fibre is a Weierstrass point. We can list all the possible 
 non-gap sequences of $p$ for any case of the list of (3) (Proposition \ref{prop: ram}).
\end{enumerate}
\end{theorem}

For the point (4) we make a crucial use of the results of Coppens \cite{cop1, cop2, cop3} and Brundu-Sacchiero \cite{BruSa} on possible gap sequences of points of trigonal curves. 
Moreover, we use the fact that the Maroni invariant of the general fibres is special $=2$ (see \cite{BPZ}).

Eventually, we find  with the same techniques, some conditions for the case $g=5$ and $u_f=4$, listed in Proposition \ref{prop: 5}.

\subsection*{Acknowledgements}
I am grateful to Pietro Pirola for many useful explanations on the subject.
I heartily thank  the referee for substantially improving the paper, both in its mathematical content (also spotting a wrong statement in the first version) and in its readability. I am partially supported by the PRIN project 20228JRCYB ``Moduli spaces and special varieties'' and by GNSAGA - INdAM. I am  a member of UMI.

\section{Fibred surfaces and their geography}\label{sec: preliminaries}
We work over the complex numbers $\mathbb C$. We will use the symbol $\cL\cong \cM$ for linear equivalence of line bundles, $L\sim M$ for linear equivalence of divisors, and $L\equiv M$ for numerical equivalence of divisors.

We call \textit{fibred surface} or sometimes simply {\em fibration} the data of a surjective morphism with connected fibres $f\colon S\rightarrow B$ from a smooth projective surface $S$ to a smooth projective curve $B$.
We denote with $b=g(B)$  the genus of the base curve. A general fibre $F$ is a smooth curve and its genus $g=g(F)$ is by definition the genus of the fibration. We consider fibrations of genus $g\geq 2$.

Some natural invariants of a fibration are the self-intersection of the relative canonical divisor $K_f=K_S-f^*K_B$, and the relative Euler characteristic $\chi_f=\chi(\cO_S)-(b-1)(g-1)$.

We say that $f$ is \textit{relatively minimal} if it does not contain any $(-1)$-curves in its fibres. 
This condition is equivalent to $K_{f}$ being a nef divisor. 
We say that a fibred surface is: \textit{smooth}  if every fibre is smooth;
 \textit{isotrivial} if all smooth fibres are mutually isomorphic;
 \textit{locally trivial} if $f$ is smooth and isotrivial;
 {\em trivial} if $S$ is birationally equivalent to $F\times B$ and $f$ corresponds to the projection on $B$. 
 If $f$ is relatively minimal this is equivalent to $S=F\times B$.
Eventually, we say that $f$ is a {\em Kodaira fibration} if it is smooth and non-isotrivial.

Let $\omega_f:=\cO_S(K_f)\cong \omega_S\otimes f^*\omega_B^{-1}$ be the relative dualizing sheaf.
The rank $g$ vector bundle $f_{*}\omega_{f}$ is called the \textit{Hodge bundle} of the fibred surface.
 By using Grothendieck-Riemann-Roch theorem we see that 
$\deg f_{*}\omega_{f}=\chi_f$.
Let us recall some important results on the Hodge bundle.

\smallskip

\begin{theorem}\label{teoremone2}
The Hodge bundle of a fibration $f\colon S\to B$ of genus $g\geq 1$  can be decomposed in two ways as a direct sum of vector sub-bundles as follows:
\begin{itemize}
\item (First Fujita decomposition \cite{Fuj1})
\begin{equation}\label{first fujita}
f_{*}\omega_{f}=\mathcal{O}_{B}^{\oplus q_{f}}\oplus \mathcal{G}, 
\end{equation}
where $\mathcal{G}$ is nef and  $H^0(B,\mathcal{G}^\vee)=H^1(B, \cG\otimes \omega_B)=0$;
\item  (Second Fujita decomposition \cite{Fuj2} \cite{CD2})
\begin{equation}\label{second fujita} 
f_{*}\omega_{f}=\mathcal{A}\oplus \mathcal{U},
\end{equation}
where $\mathcal{A}$ is ample and $\mathcal{U}$ is unitary flat.
\end{itemize}
\end{theorem}
\begin{rem}\label{rem: ossqf}
So, we see that $q_f$ is the rank of the biggest trivial subbundle of $\cE$, and $\cU$ is the biggest degree $0$ sub-bundle of $f_*\omega_f$;  Comparing the two decompositions, since every trivial bundle is unitary flat, we have that:
$$\mathcal{O}_{B}^{\oplus q_{f}}\subseteq \mathcal{U}, $$
and then it holds that $q_{f}\leq u_{f}$. 
Moreover, we see that $\chi_f=\deg \cG=\deg \cA$. As for the unitary rank, we have from the second decomposition that $u_f\leq g$, and $u_f=g$ if and only if $\chi_f=0$. If $f$ is relatively minimal, then $\chi_f=0$ is equivalent to $f$ being locally trivial.
\end{rem}

Recall two general geographical inequalities, where the left hand side is due to due to Xiao \cite{Xiao} and Cornalba-Harris \cite{C-H}, and the right hand side is classical.

\smallskip

\begin{theorem}[Cornalba-Harris, Xiao]\label{teo: geo}
For any locally non-trivial relatively minimal fibred surface of genus $g\geq 2$, we have 
\begin{equation}\label{eq: 12} 
\frac{4(g-1)}{g}\chi_f\leq K_f^2\leq 12\chi_f.
\end{equation}
Moreover, if equality holds in the left hand side, then $F$ is hyperelliptic and $q_f=0$; equality holds in the right-hand side if and only if $f$ is a Kodaira fibration.
\end{theorem}

\smallskip

Assume $f$ is relatively minimal. As mentioned above,  $q_f\leq g$, and equality holds if and only if $f$ is trivial. So, it is natural to ask: if $f$ is not trivial, what is the optimal  inequality between $g$ and $q_f$?
From now on we also assume that $f$ is non-trivial. 
Xiao initiated the study in \cite{xiao-irregular} by proving that if $b=0$ 
\begin{equation}\label{eq: xiao bound}
q_f\leq \frac{g+1}{2}. 
\end{equation}
For arbitrary $b$ he proved the bound 
\begin{equation}\label{eq: xiaobrutto}
q_f\leq \frac{5g+1}{6}.
\end{equation}
He also conjectured  that the bound  \eqref{eq: xiao bound} holds for any non-trivial fibration. 
Serrano \cite{serrano} proved that if  $f$ is isotrivial  (but not trivial), then \eqref{eq: xiao bound} holds.
Cai in  \cite{Cai} proved that if $f$ is non-isotrivial and $F$ is either hyperelliptic or bielliptic, the same bound holds. 
In \cite{konno}, Konno improved the bound \eqref{eq: xiaobrutto} to: 
\begin{equation}\label{eq: konno}
q_f \leq g \frac{5g-2}{6g-3}
\end{equation}
Pirola gave a genus $4$  counterexample to Xiao's conjecture in \cite{pirola}, by constructing for non-isotrivial fibrations a higher Abel-Jacobi map. More counterexamples have been found later by Albano and Pirola in \cite{AP}.
In  \cite{BGN}  Barja, Gonz\'alez and Naranjo proved the following: if $f$ is non-isotrivial, then $q_f\leq g-c_f$, where $c_f$ is the Clifford index of $F$.
Favale-Naranjo-Pirola \cite{FNP} proved the stronger inequality $q_f\leq g-c_f-1$ for families of plane curves of degree $d \geq 5$.

It is natural to ask wether these bounds can be extended to the unitary rank $u_f$. 
In the paper \cite{GTS} the author with Gonz\'alez and Torelli extended the results of  \cite{BGN} and  \cite{FNP} to the unitary rank: in particular we proved
the inequality $u_f \leq g - c_f.$
Moreover, if the general fibre is a plane curve
of degree $d \geq 5$, then
$ u_f \leq  g - c_f - 1.$

In \cite{Sto24}  we proved the following generalization of Konno's bound: 
\begin{equation}\label{eq: konno-uf}
u_f < g \frac{5g-2}{6g-3}
\end{equation}

\begin{rem}\label{rem: confronto}
If $\cU$ has finite monodromy, then up to an \'etale base change it becomes trivial, and so the unitary rank becomes equal to the relative irregularity. 
However, the unitary summand $\cU$ of the Hodge bundle can have infinite monodromy as proved by Catanese and Dettweiler in \cite{CD2}. In these cases the relative irregularity is strictly smaller than the unitary rank  also up to base changes.
 \end{rem}

\subsection{Xiao's set-up, and a slope inequality}
Let us briefly recall that the Harder-Narasimhan sequence of  a vector bundle $\cE$ over a curve $B$:
it is the unique filtration of subbundles
$$0=\cE_0\subset \cE_1\subset \ldots \subset \cE_l=\cE$$
satisfying the following assumptions:
\begin{itemize}
\item  $\cE_i/\cE_{i-1}$ is  a $\mu$-semistable vector bundle for any $i=1, \ldots l$;
\item if we set $\mu_i:= \mu(\cE_i/\cE_{i-1})$ for any $i=1, \ldots l$, we have that $\mu_i<\mu_{i-1}$.
\end{itemize}
Note that $\mu_1>\mu(\cE)>\mu_l$, unless $\cE$ is $\mu$-semistable, in which case $1=l$ and $\mu(\cE)=\mu_1$. 


Let us consider the case of a fibred surface $f\colon S\to B$, and consider its Hodge bundle  $\cE=f_*\omega_f$.
As $\cE$ is nef by Theorem \ref{teoremone2}, we have that $\mu_l$ is greater or equal to $0$. 

\smallskip

\begin{rem}\label{rem: ok}
The unitary flat summand $\cU$ in \eqref{second fujita} is the biggest subbundle of $\cE$ with degree $0$. 
So, in particular we have that $\mu_{l}=0$ if and only if $\mathcal U\not=0$ if and only if $u_f\not =0$.
Moreover, if $\mu_l=0$, then $\cE_{l-1}$ is precisely $\mathcal A$, if $\mu_l>0$, then $\mathcal A$ is the whole Hodge bundle.
\end{rem}

We follow the  set-up given by Xiao in his seminal paper  \cite{Xiao}. 
Let us consider the Hodge bundle $\cE=f_*\omega_f$, and its Harder-Narasimhan sequence as described above.
For any $i=1,\ldots l$, the sheaf homomorphism $\eta_i$
$$
\xymatrix{f^*\cE_i\ar@/^1pc/[rr]^{\eta_i} \ar[r]& f^*f_*\omega_f\ar[r]&\omega_f}
$$
induces a rational map $\phi_i\colon S\dasharrow \pr_B(\cE_i)$. 
Let $\cO_{\pr}(1)$ be  the tautological line bundle of $\pr_B(\cE_i)$, and let $H_i$ be any associated divisor.  
Let $M_i:=\phi_i^*(H_i)$ on $S$.
%
Moreover, for any $i=1,\ldots l$, call $Z_i$ the divisorial base locus of $\phi_i$, i.e. the effective divisor in $S$ such that  
$\cE_i\subseteq  f_*\omega_f(-Z_i)$ and such that 
the evaluation homomorphism
\[f^*\cE_i\longrightarrow f^*f_*\omega_f(-Z_i)\longrightarrow \omega_f(-Z_i)\]
is surjective in codimension $2$. 
Of course
\[Z_1\geq Z_2\geq \ldots \geq Z_l=0.\]
\begin{definition}
With the above notations, for any $i=1,\ldots l$, we denote by $a_i:=Z_i\cdot F$; in other words, $a_i$ is the degree of the base locus of the linear system $\cE_i\otimes\kapp(t)$ as a sub-system of $\cE\otimes\kapp(t)=H^0(F_t,\omega_{F_t})$ for $t\in B$ general.
Let moreover $Z_i=\sum_j m_jG_j$ be the decomposition of $Z_i$ in  irreducible components. 
We define 
\[\alpha_i:=\max_j\{m_j\mid  F\cdot G_j>0\}.\]
\end{definition}
\begin{rem}
Note that we have 
\[a_1\geq a_2\geq\ldots \geq a_l=0, \quad \alpha_1\geq \alpha_2\geq\ldots \geq \alpha_l=0. \]
 Note moreover that for any $i=1,\ldots l$, we have $\alpha_i\leq a_i$, and if in some cases we have equality then $Z_i=\alpha_i\Gamma$ where $\Gamma$ is a section of $f$.
\end{rem}
In \cite{Sto24}, we proved the following:

\smallskip

\begin{theorem}\label{lem: konnoplus}
Let $f\colon S\to B$ a relatively minimal non-locally trivial fibred surface of genus $g\geq 2$. Let $a:=a_1$, and $\alpha:=\alpha_1$.
The following inequality holds:
\begin{equation}\label{eq: 2}
K_f^2\geq \left(\frac{(\alpha+1)(4g-4-a)+a}{\alpha+1}\right)\mu_1\geq \left(\frac{(\alpha+1)(4g-4-a)+a}{\alpha+1}\right)\frac{\chi_f}{g-u_f}.
\end{equation}
\end{theorem}

\section{Main results}\label{sec: 6}

Let us first collect here some results on the extremal case, mainly due to Pirola and his collaborators:

\smallskip

 \begin{theorem}[Beorchia, Pirola, Serrano, Stoppino, Zucconi]\label{teo: estremale}
Let $f\colon S\to B$ a relatively minimal non trivial fibred surface of genus $g\geq 4$ such that $u_f= g-1$. 
Then the following hold:
\begin{itemize}
\item[]
\begin{enumerate}
\item [(1)] (\cite[Theorem 3.6 and Remark 3.7]{pietroBUMI}) up to an \`etale base change $B'\to B$ the induced fibration $f'$ has $q_{f'}=u_{f'}$; 
\item [(2)] (\cite[Thm 5.2]{Sto24}) the smooth fibres $F_t$ are a covering of a non-isotrivial family of elliptic curves;
\item [(3)] (\cite[Thm 5.2]{Sto24}) the fibration is not Kodaira , i.e. $K_f^2<12\chi_f$;  
\item [(4)] (\cite{BPZ}) the general fibre $F$ is trigonal and has special Maroni invariant. 
\end{enumerate}
\end{itemize}
\end{theorem}
In  \cite[Section 10]{CFP} some other explicit conditions are found on the case $g=6$.

Let $f\colon S\to B$ be a relatively minimal locally non-trivial fibred surface of genus $g$ and $u_f= g-1$. 
Then the Harder-Narasimhan sequence of the Hodge bundle $\cE=f_*\omega_f$ is 
\[0=\cE_0\subset \cE_1=\cA\subset \cE_2=\cE,\]
and $r_1=\rk\cE_1=1$, so $\cA$ is an ample line bundle on $B$.
 Let us use in this case the notation  $Z:=Z_1$ and $\alpha:=\alpha_1$.
 
 \smallskip
 
\begin{lemma}\label{lem: serve}
Under the assumptions above, we have $\omega_f\cong f^*\cA(Z)$.
\end{lemma}
\begin{proof}
The map $\phi_1\colon S \to \pr_B(\cE_1)\cong B$ is a morphism in this case, and it coincides with the fibration $f$. 
The tautological line bundle on $\pr_B(\cE_1)$ is  $\cE_1=\cA$.
The evaluation morphism 
\[
f^*\cA\longrightarrow \omega_f(-Z)
\]
is a surjective morphism between line bundles, so it is an isomorphism, as wanted.
\end{proof}
We use Theorem \ref{lem: konnoplus} to give two necessary conditions for the case $g=6$.

\smallskip

\begin{proposition}\label{prop: 6}
Let $f\colon S\to B$ a relatively minimal locally non-trivial fibred surface of genus $g= 6$ such that $u_f= 5$. 
Then  $S$ has positive  index $\tau(S)>0$
and  $\alpha\geq 5$.
\end{proposition}
\begin{proof}
Using Theorem \ref{lem: konnoplus}, we obtain 
\[K_f^2\geq \left( 2g-2+\frac{2g-2}{\alpha+1}\right)\mu_1\geq \left(10+\frac{10}{\alpha+1}\right)\chi_f.\]
This implies that the index $\tau(S)=K_S^2-8\chi(\cO_S)=K_f^2-8\chi_f$  necessarily is positive. 

In order to prove the second statement, observe that if $\alpha\leq 3$, the inequality above would contradict  \eqref{eq: 12}.
If $\alpha=4$, we would have 
$K_f^2\geq 12\chi_f$, so necessarily equality holds and $f$ is a Kodaira fibration,  but this contradicts point (3) of  Theorem \ref{teo: estremale} 
\end{proof}
\begin{rem}
So, we have obtained that under the assumptions of Proposition \ref{prop: 6}, $Z$ necessarily has an horizontal component $\Gamma$ of multiplicity greater or equal to $5$. 
If $\alpha\geq 6$ then  this component is a section of $f$. 
Indeed, calling $b:= \Gamma\cdot F$, we have 
\[
2g-2=10=K_f\cdot F\geq Z\cdot F\geq \alpha b.
\]
Hence, if $\alpha\geq 6$ then $b=1$, i.e. $\Gamma$ is a section.
In the case $\alpha=5$ we could have that $Z=5\Sigma$, where $\Sigma$ is a $2:1$ covering of $B$. Unfortunately we have not been able to exclude this case.
\end{rem}

We now see that in the extremal case $u_f =g-1$, we can improve quite strongly the slope inequality of Theorem \ref{lem: konnoplus} via a direct computation as follows:

\smallskip

\begin{proposition}\label{prop: bound estremale}
Let $f\colon S\to B$ a relatively minimal  locally non-trivial fibred surface of genus $g\geq 2$ such that $u_f= g-1$. Let $Z=\sum^k_{j=1}m_jG_j+Z'$, where the $G_j$'s are the irreducible components of $Z$ which are not vertical with respect to $f$ and the $m_j$'s are their respective multiplicities. 
Let $b_j:=G_j\cdot F$. 
Then we have the following inequalities:
\begin{equation}
K_f\cdot Z\geq \left( \sum_{j=1}^k\frac{m_j}{m_j+1}b_j\right)\mu_1.
\end{equation}
\begin{equation}\label{eq: meglio}
K_f^2\geq \left(2(g-1)+ \sum_{j=1}^k\frac{m_j}{m_j+1}b_j\right)\mu_1.
\end{equation}
\end{proposition}
\begin{proof}
Using Theorem \ref{teo: estremale} we can suppose that $q_f=u_f$, because the ratio $K_f^2/\chi_f$ remains the same under \'etale base changes.
So, under our assumptions, $f_*\omega_f =\cA\oplus{\cO_B}^{\oplus (g-1)}$.

Observe that $\deg f_*\cO_S(Z)=-(g-1)\mu_1$.
Indeed, using   Lemma \ref{lem: serve}
\[f_*\omega_f\cong f_*(\cO_S(Z))\otimes f^*\cA)\cong f_*\cO_S(Z)\otimes \cA,\]
so the formula for the degree of vector bundles gives: 
\[\deg f_*\cO_S(Z)= -g\deg\cA +\deg f_*\omega_f=-g\chi_f +\chi_f=-(g-1)\mu_1,\] 
as wanted.
Now observe that by the relative Serre duality we have 
\[R^1f_*(\cO_S(Z))\cong (R^1f_*\omega_f)\otimes \cA^{-1}\cong \cO_S\otimes \cA^{-1}\cong \cA^{-1},\]
 so $\deg R^1 f_*(\cO_S(Z))=-\mu_1$.
Let us use Grothendieck-Riemann-Roch
 for $\cO_S(Z)$:
 \[
 \deg f_*\cO_S(Z)=\frac{{Z}^2}{2}-\frac{K_f\cdot Z}{2}+\chi_f+\deg R^1 f_*(\cO_S(Z))=\frac{{Z}^2}{2}-\frac{K_f\cdot Z}{2}.
 \]
Summing up, we have 
$ {Z}^2=K_f\cdot Z-2(g-1)\mu_1.$
Now we compute $K_f^2$ using the equality above and fact that $K_f\cong Z+f^*\cA\equiv Z+\mu_1F$:
 \begin{equation}\label{eq: equa}
 K_f^2={Z}^2+4(g-1)\mu_1= K_f\cdot Z+ 2(g-1)\mu_1.
 \end{equation}
 Next, we want to estimate the number $K_f\cdot G_i$ for any $i=1,\ldots , k$.
 By making explicit the terms, we have 
 \[K_f\cdot G_i=\left(\sum^k_{j=1}m_jG_j+Z'+\mu_1F\right)\cdot G_i=\sum^k_{j=1}m_j(G_j\cdot G_i)+(Z'\cdot G_i)+\mu_1(F\cdot G_i)\geq m_iG_i^2+ b_i\mu_1.\]
Now, observe that $(K_f+G_i)\cdot G_i=\deg r(\pi)$, where $\pi\colon G_i\to B$ is the morphism induced by restriction of $f$, and $r(\pi)$ is its ramification divisor. 
So, $K_f\cdot G_i\geq -G_i^2$.
Combining this inequality with the above formula, we obtain
\[
(m_i+1)K_f\cdot G_i\geq b_i\mu_1.
\]
Now compute $K_f^2$ using \eqref{eq: equa}:
\[
K_f^2\geq 2(g-1)\mu_1+K_f\cdot Z=2(g-1)\mu_1+K_f\cdot \left(\sum^k_{j=1}m_jG_j+Z'\right)\geq 2(g-1)\mu_1+\left(\sum^k_{j=1}\frac{m_jb_j}{m_j+1}\right)\mu_1,
\]
as wanted.
\end{proof}

\begin{rem}\label{rem: migliore}
Although the function $m\mapsto \frac{m}{m+1}$ is increasing with $m$, the contribution $ \sum_{j=1}^k\frac{m_j}{m_j+1}b_j$ appearing in inequality \eqref{eq: meglio} is better if the multiplicities of the horizontal irreducible components of $Z$ are low, as we now illustrate.
Recall that $\sum_{j=1}^{k}b_jm_j=2g-2=Z\cdot F$. 
Let us consider the case when $m_j=1$ for all $j=1,\ldots, k$. Then $\sum_{j=1}^kb_j=2g-2$ and we  have 
\[
\sum_{j=1}^k\frac{m_j}{m_j+1}b_j=\frac{1}{2}\sum_{j=1}^kb_j=\frac{2g-2}{2}=g-1.
\]
If $m_1\geq2$ and $m_j=1$ for all $j=2,\ldots, k$, then 
\begin{align*}
 \sum_{j=1}^k\frac{m_j}{m_j+1}b_j=\frac{m_1b_1}{m_1+1}+ \sum_{j=2}^k\frac{m_j}{m_j+1}b_j=\frac{m_1b_1}{m_1+1}+\frac{2g-2-b_1m_1}{2}=\\
=\frac{2g-2}{2} -b_1m_1\left(\frac{1}{2}-\frac{1}{m_1+1}\right)<g-1.
 \end{align*}
 So, it is clear that the maximal possible contribution is $g-1$ if and only if  $Z$ is reduced, and that the bigger the $m_j$'s are, the smaller the contribution is.
 \end{rem}

 \smallskip 
 
Now we use the strong inequality above in order to make a skimming of the possible classes of $K_f$. Note that we use here numerical equivalence because $f^*\cA\equiv \mu_1F$ but it is not necessarily linearly equivalent. However, if we replace $\mu_1F$ with $f^*\cA$, we have linear equivalence in the statement below.

\smallskip

\begin{proposition}\label{prop: elenco}
Let $f\colon S\to B$ a relatively minimal locally non-trivial fibred surface of genus $g= 6$ such that $u_f= 5$. 
The possibilities for the class of $K_f$ are the following. 
\begin{itemize}
\item[]
\begin{itemize}
\item[(1)] $\alpha=10$ and $K_f\equiv 10\Gamma +Z'+\mu_1F$ where $\Gamma$ is a section of $f$.
\item[(2)] $\alpha=9$ and $K_f\equiv 9\Gamma +\Gamma' +Z'+\mu_1F$ where $\Gamma$ and $\Gamma'$ are distinct sections of $f$.
\item[(3)] $\alpha=8$ and 
   \begin{itemize}
    \item either $K_f\equiv 8\Gamma +2\Gamma'+Z'+\mu_1F$, where $\Gamma$ and $\Gamma' $  are distinct sections of $f$.
    \item or $K_f\equiv 8\Gamma +\Gamma'+\Gamma''+Z'+\mu_1F$, where $\Gamma, \Gamma' $ and $\Gamma''$ are distinct sections of $f$.
    \item or  $K_f\equiv 8\Gamma +\Sigma +Z'+\mu_1F$ where  $\Sigma$ a bisection of $f$.
     \end{itemize}
\item[(4)] $\alpha=7$ and  $K_f\equiv 7 \Gamma + 3\Gamma'+Z'+\mu_1F $, where $\Gamma, \Gamma'$ are  distinct sections of $f$,
\item[(5)] $\alpha =6$ and 
$K_f\equiv 6\Gamma + 4\Gamma'+Z'+\mu_1F $, where $\Gamma, \Gamma'$ are  distinct sections of $f$. 
\item[(6)] $\alpha=5$ and 
\begin{itemize}
\item  either $K_f\equiv 5\Gamma + 5\Gamma'+Z'+\mu_1F $, where $\Gamma, \Gamma'$ are distinct sections of $f$,
\item or $K_f\equiv 5\Sigma+Z' +\mu_1F$, where $\Sigma $ is a bisection of $f$.
\end{itemize}
\end{itemize}
In all the above $Z'$ is the vertical part of $Z$, if any.
\end{itemize}
\end{proposition}
\begin{proof}
Let $G:=\sum_{j=1}^{k}G_j$ be the union of the irreducible horizontal components of $Z$. 
 We claim that if $G\cdot F\geq 4$, then we get a contradiction. 
Indeed, from \eqref{eq: meglio} we would obtain a contribution to the coefficient of $\mu_1$ 
of the form 
\[\sum_{j=1}^k\frac{m_j}{m_j+1}b_j\geq \frac{1}{2}G\cdot F\geq 2,\]
and so \eqref{eq: meglio} would give $K_f^2\geq 12\chi_f $, and obtain a contradiction using point (4) of 
Theorem \ref{teo: estremale}.

The remaining possible cases according to the values of $\alpha$ are the following. In the list below  $\Gamma,\Gamma'$ and $\Gamma''$ are  distinct sections of $f$, $\Sigma$ is a multisection, i.e. a degree $d>1$ covering  of the base $B$. The degree is specified case by case. We omit here the vertical part $Z'$ for simplicity, because it is not involved in the formula of Theorem \ref{teo: estremale}.
\begin{itemize}
\item[]
\begin{itemize}
\item[$(1)'$] $\alpha=10$ and $K_f\equiv 10\Gamma +\mu_1F$,
\item[$(2)'$] $\alpha=9$ and $K_f\equiv 9\Gamma +\Gamma' +\mu_1F$,
\item[$(3)'$] $\alpha=8$ and 
   \begin{itemize}
    \item[$(3a)'$] $K_f\equiv 8\Gamma +2\Gamma'+\mu_1F$, 
    \item[$(3b)'$]  $K_f\equiv 8\Gamma +\Gamma'+\Gamma''+\mu_1F$,
    \item[$(3c)'$]  $K_f\equiv 8\Gamma +\Sigma +\mu_1F$ where  $\Sigma$ a bisection of $f$.
     \end{itemize}
\item[$(4)'$] $\alpha=7$ and 
 \begin{itemize}
    \item[$(4a)'$]  $K_f\equiv 7\Gamma +\Gamma'+2\Gamma''+\mu_1F$, 
    \item[$(4b)'$]  $K_f\equiv 7\Gamma +3\Gamma'+\mu_1F$, 
    \end{itemize}
\item[$(5)'$] $\alpha =6$ and 
 \begin{itemize}
      \item[$(5a)'$] $K_f\equiv 6\Gamma + 3\Gamma + \Gamma'+\mu_1F $,
     \item[$(5b)'$]  $K_f\equiv 6\Gamma + 4\Gamma'+\mu_1F $,
      \item[$(5c)'$] $K_f\equiv 6\Gamma + 2\Sigma+\mu_1F $, where $\Sigma $ is a bisection,
    \end{itemize} 
    \item[$(6)'$] $\alpha=5$ and  
    \begin{itemize} 
     \item[$(6a)'$] $K_f\equiv 5\Gamma + 2\Gamma'+3\Gamma''+\mu_1F $,
     \item[$(6b)'$] $K_f\equiv 5\Gamma + \Gamma' +4\Gamma''+\mu_1F $,
     \item[$(6c)'$] $K_f\equiv 5\Gamma + 5\Gamma'+\mu_1F $,
     \item [$(6d)'$] $K_f\equiv 5\Sigma+\mu_1F $ where $\Sigma$ is a bisection of $f$.
     \end{itemize}
\end{itemize}
\end{itemize}
 We are going to rule out the great part of these possibilities by using Theorem \ref{teo: estremale}.
Let us  consider  the cases we want to exclude:
 \begin{itemize}
 \item[]
 \begin{enumerate}
\item[$(4a)'$]  Suppose $K_f\equiv 7\Gamma +\Gamma'+2\Gamma''+\mu_1F$.
Using inequality \eqref{eq: meglio} we obtain
\[K_f^2\geq 10\mu_1+K_f\cdot Z\geq  \left(10+\frac{7}{8}+\frac{1}{2}+1\right)\mu_1>12\mu_1.\]

\item[$(5a)'$] $K_f\equiv 6\Gamma + 3\Gamma + \Gamma'+\mu_1F $ would imply 
 \[
 K_f^2=K_f\cdot Z+ 10 \mu_1\geq \left(10+\frac{6}{7}+\frac{3}{4}+\frac{1}{2}\right)\mu_1
 =\left(10+\frac{57}{28}\right)\mu_1>12\mu_1. \]

\item[$(5b)'$] $K_f\equiv 6\Gamma + 2\Sigma+\mu_1F $, where $\Sigma $ is a bisection,
 would imply 
 \[
 K_f^2=K_f\cdot Z+ 10 \mu_1\geq \left(10+\frac{6}{7}+\frac{4}{3}\right)\mu_1
>12\mu_1. \]

\item[$(6a)'$] $K_f\equiv 5\Gamma + 2\Gamma'+3\Gamma''+\mu_1F $,
we obtain 
\[
 K_f^2=K_f\cdot Z+ 10 \mu_1\geq \left(10+\frac{5}{6}+\frac{2}{3}+\frac{3}{4}\right)\mu_1>12\mu_1.
 \]
 \item[$(6b)'$] $K_f\equiv 5\Gamma + \Sigma +3\Gamma'+\mu_1F $, where $\Sigma$ is a  bisection,    then we obtain 
\[
 K_f^2=K_f\cdot Z+ 10 \mu_1\geq \left(10+\frac{5}{6}+1+\frac{3}{2}\right)\mu_1>12\mu_1.
 \]
     \item[$(6c)'$] $K_f\equiv 5\Gamma + \Gamma' +4\Gamma''+\mu_1F $,    then we obtain 
\[
 K_f^2= K_f\cdot Z+ 10 \mu_1\geq \left(10+\frac{5}{6}+\frac{1}{2}+\frac{4}{2}\right)\mu_1>12\mu_1.
 \]
   \end{enumerate} 
   \end{itemize}
   \end{proof}

So far we  have excluded many numerical classes for $K_f$ just using Proposition \ref{prop: bound estremale} combined with the fact that the slope $K_f^2/\chi_f$ has to be strictly smaller than $12$. Now we make the trigonal structure of $F$ come into play.
We see that in the cases when $\alpha\geq 6$  we have a control on the possible Weierstrass sequence of the intersection of the section $\Gamma$ with the general fibre $F$.
We first need to recall this classical result.

\smallskip

\begin{lemma}\label{lem: maroni}
Let $F$ be a trigonal genus $6$ curve with and let calling $D$ a divisor belonging to the unique $g^1_3$.
The following hold:
\begin{itemize}
\item[]
\begin{itemize}
\item[(i)] $F$ has general Maroni invariant $0$ if and only if $K_f\sim 2D+R$, where $|R|$ is a (unique) base point-free $g^1_4$.
\item[(ii)] $F$ has special Maroni invariant $2$ if and only if  $K_F\sim  3D+q$, for some $q\in F$.
\end{itemize}
\end{itemize}
Moreover, $F$ has general Maroni invariant if and only if it has a base-point-free $g^1_4$.
\begin{proof}
(i) $F$ has Maroni invariant $0$ if and only if its  canonical image lies in the Hirzebruch surface  $\F_0=\pr^1\times \pr^1$ as a curve of bi-degree $(3,4)$. Call $L_1$ the fibre of the first ruling and $L_2$ the fibre of the second one. The  first ruling defines the $g^1_3=|D|$ and the second ruling defines a base-point-free $g^1_4=|R|$ over $F$. In other words, ${L_1}_{|F}\sim R$ and ${L_2}_{|F}\sim D$.
By adjunction 
\[
K_F\sim (K_{\pr^1\times \pr^1}+F)_{|F}\sim(L_1+2L_2)_{|F}\sim 2D+R,
\]
as claimed. On the other hand, if a genus $6$ trigonal curve has a base point free $g^1_4$, then we can define a map 
$\varphi \colon F\longrightarrow \pr^1\times \pr^1$ combining the two pencils. As the $g^1_4$ is base point free, it is not composed with the $g^1_3$ and so the morphism $\varphi$ is finite. By the genus formula it is immediate to see that it is an embedding.
 
\noindent (ii) $F$ has Maroni invariant $2$ if and only if the canonical image of $F$ lies in a Hirzebruch surface  $\F_2$. 
Recall that $\F_2$ is fibred over $\pr^1$, and call $\pi$ this fibration.
Let us call $\Sigma$ the negative section of $\F_2$ and $K$ a fibre of $\pi$. So we have 
\[\Sigma^2=-2, \,\,\, \Sigma\cdot K=1,\,\,\, K^2=0.
\]
The curve $F$ belongs to the linear system $\mid 3\Sigma+ 7 K\mid$, as can be easily checked. 
The $g^1_3=|D|$ is cut out on $F$ by the fibres of $\pi$; in other words $K_{|F}\sim D$. 
By the adjunction theorem $K_F\sim(\Sigma+3K)_{|F}\sim 3D+q$, where $q=\Sigma\cap F$, as wanted. 
On the other hand, if $F$ is trigonal of genus $6$ and $K_F\sim  3D+q$, for some $q\in F$, then necessarily $F$ has Maroni invariant 2. 
Indeed, if  $F$ would have Maroni invariant $0$, we would have  $ 3D+q\sim 2D+R$, equivalently $D+q\sim R$, a contradiction because $R$ is base-point-free.

The last statement follows from (i).
\end{proof}
\end{lemma}
\begin{proposition}\label{prop: ram}
In the situation of Proposition \ref{prop: elenco}, if $\alpha>5$ then  the point $p:=\Gamma\cap F$ is a Weierstrass point for the curve $F$. 
We have the following possibilities:
\begin{itemize}
\item[]
\begin{itemize}
\item[(1)]  $p$ is a total ramification point for the trigonal morphism of  $F$, whose non-gap sequence is:
$(1,\,2,\,4,\,5,\,7,\,10)$;  this can only happen for $\alpha=9 $.
\item[(3)]  $p$ is not a ramification point of the trigonal morphism; in this case the possible  non-gap sequences of $p$ are:
\begin{enumerate}
\item $(1,\,2,\,3,\, 4,\, 5,\,11)$; this can only happen for $\alpha= 10$.
\item $(1,\,2,\,3,\, 4,\, 5,\,9)$; this can only happen for $\alpha= 8$.
 \item $(1,\,2,\,3,\, 4,\, 5,\,8)$; this can only happen for $\alpha=7$.
 \item $ (1,\,2,\,3,\, 4,\, 5,\,7)$  this can only happen for $\alpha=6$.
\end{enumerate}
\end{itemize}
\end{itemize}
\end{proposition}
\begin{proof}
First we prove that $p$ is a Weierstrass point of $F$ if $\alpha\geq 6$: indeed in this case we see that $h^0(6p)=h^0(K_F-R)$, where $R$ is an effective divisor on $F$ of degree $4$, so it holds $h^0(6p)\geq h^0(K_F)-4=2$.

We want to determine the possible non-gap sequences of the point $p$.
Recall that we know that the fibre $F$ is trigonal with special Maroni invariant $2$ \cite{BPZ}. 
Thanks to the work of Coppens \cite{cop1, cop2, cop3} and Brundu-Sacchiero \cite{BruSa} we have a complete description of the possible non-gap sequences of a Weierstrass point on a trigonal curve. 
We list here the possibilities. 
\begin{enumerate}
\item \cite{cop1, cop3} if $p$ is a total ramification point of the trigonal morphism (i.e. $|3p|=g^1_3$), then the possible non-gap sequence are:
\begin{equation}\label{questa}
(1,\,2,\,4,\, 5,\, 7,\,10)
\end{equation}
\begin{equation}\label{quella}
(1,\,2,\,4,\, 5,\, 8,\,11)
\end{equation}
\item \cite[Theorem 3.1]{BruSa} \cite{cop3} If $p$ is a simple ramification point of of the trigonal morphism (i.e. $|2p+r|=g^1_3$, with $r\not=p$), then the possible non-gap sequence in $p$ are:
\begin{equation}
(1,\,2,\,3,\, 4,\, 5,\,7)
\end{equation}
\begin{equation}\label{ram}
(1,\,2,\,3,\, 4,\, 6,\,8)
\end{equation}
\item \cite[Theorem 3.5]{BruSa} If $p$ is not a ramification point of of the trigonal morphism  then there exist integers $r,\gamma$, with $4\leq r\leq 5$,  
$1\leq \gamma \leq 2r-5$  such that the possible non-gap sequence in $p$ splits into two sequences of consecutive integers as follows:
\begin{equation}\label{eq: laltra}
(1,\,2,\,\ldots,r, r+1+\gamma,\ldots ,6+\gamma)
\end{equation}
\end{enumerate}
Let us start with case $\alpha=10$. In this case $p$ is a subcanonical point of $F$, i.e. $K_F\sim(2g-2)p=10p$. 
So we have 
\[h^0(11p)=h^0(K_F(p))=h^0(K_F)= h^0(10p)=6,\]
hence the last non-gap is $11$. 
So the only possibilities for the sequences are \eqref{quella} or \eqref{eq: laltra} with $r=\gamma=5$. 
The last sequence is 
\[(1,\,2,\,3,\, 4,\, 5,\,11).\]
The possible sequences for subcanonical points on curves of genus $g\leq 6$ have been classified by Bullock in the table of \cite[Sec.4]{bull}. From this result, we see that the sequence \eqref{quella} can not happen.

Let us now consider case $\alpha= 9$. 
So $K_F\sim9p+ q$, where $p\not = q$ (because the sections $\Gamma $ and $\Gamma'$ are distinct). 
Observe that $h^0(11p)=6$ by Riemann-Roch and  $  h^0(10p)=h^0(9p)=5$, because $10p\not\sim K_F$, and $9p\sim K_F-q$. 
So, the last non-gap number is $10$. 
We see immediately from the lists above that $p$ cannot be a simple ramification point. 
The only possibilities are the sequences \eqref{questa}, or \eqref{eq: laltra} with $r=5, \,\gamma=4$. 
In this last case the sequence is 
\[(1,\,2,\,3,\, 4,\, 5,\,10).\]


Consider the case $\alpha=8$. So, $K_F\sim 8p+ R$, where $R$ is effective of degree 2 not containig $p$ in its support.  
So, $h^0(9p)=4+h^0(R-p)=4$, and $h^0(8p)=h^0(K_F-R)=4$, so the last non-gap  necessarily is $9$. 
This implies that $p$ cannot be of total nor of simple ramification, and the only possibility is 
\eqref{eq: laltra} with $\gamma=3$ and $r=4$ or $r=5$. 
These sequences are:
\[(1,\,2,\,3,\, 4,\, 8,\,9)\,\,
\mbox{ and }\,\,\,\,\, (1,\,2,\,3,\, 4,\, 5,\,9).\]
We see that the first sequence is impossible. Assume by contradiction that it holds. We would have $4=h^0(7p)=2+ h^0(p+R)$, so $|p+R|=g^1_3$.
In this case the divisor $7p$ gives an generically injective morphism  $\varphi\colon F\to \pr^3$. 
Indeed, $h^0(7p-r-r')=3$ if and only if $h^0(p+R+r+r')=3$ so if and only if $F$ has a $g^2_5$, but this implies that $F$ is a plane quintic, so it has gonality $4$, a contradiction.
By Castelnuovo's bound (see for instance \cite{Harris}) for irreducible non-degenerate curves, the arithmetic genus of the image curve $\overline F$ is smaller or equal to $6$. 
So the genus is precisely $6$, $\varphi$ is an isomorphism on its image, and $F$ is a  curve reaching Castelnuovo's bound. 
Then, again by the result of Castelnuovo,  $F$ is contained in a quadric in $\pr^3$, so in $\pr^1\times \pr^1$, and so it has Maroni invariant $0$ by Lemma \ref{lem: maroni}, a contradiction.

Let us consider the case $\alpha=7$. 
By Proposition  \ref{prop: elenco} we have  $K_F\sim 7p+ 3q$, where $p\not =q$.
Let us first consider the cases where $p$ is  a ramification point.
If $p$ is of total ramification, then by Lemma \ref{lem: maroni} we have $7p+3q\sim 9p+t$ for some $t\in F$. This implies 
that $3q\sim 2p+t$ so the curve is either hyperelliptic if $t=q$, or $|3q|=|2p+t|=g^1_3$ and so $p$ is of simple ramification. 
Both are contradictions. 
If $p$ is of simple ramification, we would have $7p+3q\sim 6p+3r+t$ for some $r,t\in F$ such that $|2p+r|=g^1_3$.
So, we have $p+3q\sim 3r+t$ a $g^1_4$ that  by Lemma \ref{lem: maroni} needs to have a base point. 
Recall that  $r\not=p$
If $p=t$ and $q\not =r$ then we would have $|3r|=g^1_3=|2p+r|$, absurd.
On the other hand, if  $q=r$ this implies $p\sim t$ and so $p=t$, but then $p$ would be of total ramification, a contradiction. 
As for the case when $p$ is  non-ramified, it is not hard to see that 
the only possibilities are
\eqref{eq: laltra} with $\gamma=2$ and  $r=5$ or or  $\gamma=3$ and $r=4$:
\[(1,\,2,\,3,\, 4,\, 5,\,8)\,\,\,\mbox{ and }\,\,\,\,\, (1,\,2,\,3,\, 4,\, 8,\,9).\]
In the second sequence, we have  $h^0(5p)=2=h^0(2p+3q)$ and $3=h^0(6p)=1+ h^0(p+3q)$ which means that $|p+2q|=g^1_3$. Then by Lemma \ref{lem: maroni} we would have $7p+3q\sim 3(p+2q)+t$ for some $t\in F$, so $4p\sim 3q+t$ which gives a contradiction.

Let us now consider the case $\alpha=6$. From Proposition \ref{prop: elenco} we see that in this case $K_f\sim 6p+4q$, with $p\not =q$. 
If $p$ is of total ramification we have by Lemma \ref{lem: maroni} that $9p+t\sim K_F\sim 6p+4q$ so $3p+t\sim 4q$, so we need to have $t=q$ and $3q\sim 3p$; this  implies that $10q\sim K_F$, so $q$ is a subcanonical point  for $F$ of total ramification. This is impossible by the result of Bullock recalled in the case  $\alpha=10$.
If $p$ is of simple ramification we  have $6p+4q\sim 6p+3r+t$ for some $r,t\in F$. So, we have $4q\sim 3r+t$ and again this necessarily implies $t=q$ and $3q$ is the $g^1_3$ on $F$. But this implies that $2p+r\sim 3r$, a contradiction. 
Assume now that  $p$ is a non-ramified point, we have necessarily that $h^0(8p)=3$.
 Indeed, $h^0(4q-2p)=0$: if this was not the case, we would have a degree $2$ effective divisor $T$ such that $T\sim 4q-2p$, so $4q\sim 2p+T$. This implies by Lemma \ref{lem: maroni} that $T $ contains $q$ in its support, and either $2q\sim 2p$, which would imply that $F$ is hyperelliptic, or  $3q\sim 2p+t$, contrary to the assumption that $p$ is of simple ramification.
Now consider $h^0(7p)=2+h^0(4q-p)$. With similar arguments we can see that necessarily $h^0(4q-p)=0$, so $h^0(7p)=2=h^0(6p)$.
The only possibility here is the sequence \eqref{eq: laltra} with $\gamma=1$ and  $r=5$:
\[  (1,\,2,\,3,\, 4,\, 5,\,7) \]
\end{proof}


We end this section with a similar result on genus five maximally irregular fibrations.
\begin{proposition}\label{prop: 5}
Let $f\colon S\to B$ a relatively minimal locally non-trivial fibred surface of genus $g= 5$ such that $u_f= 4$. 
Then $f$ is a trigonal fibration, $S$ has positive  index $\tau(S)>0$, and $Z$ has at least one multiple horizontal component, i.e. $\alpha\geq 2$.
\end{proposition}
\begin{proof} 
The divisor  $Z$ necessarily has  horizontal components, because $Z\cdot F=2g-2=8$.  If all of them are reduced, ee use \eqref{eq: meglio} of Proposition \ref{prop: bound estremale} and obtain a slope greater than or equal to $12$. 
Indeed, let $Z=\sum_jG_j+Z'$, with $Z'$ vertical (and $m_j=1$). We have:
\[
K_f^2\geq \left(2(g-1)+ \sum_{j=1}^k\frac{1}{2}b_j\right)\mu_1= \left(8+\frac{8}{2}\right)\mu_1=12\chi_f,
\] 
because $\sum_jb_j=2g-2=8$.
So, we obtain a contradiction with point (3) of Theorem \ref{teo: estremale}.
\end{proof}

\end{document}